\documentclass{article}
\usepackage{amsthm}
\usepackage[english]{babel}
\usepackage[letterpaper,top=2cm,bottom=2cm,left=3cm,right=3cm,marginparwidth=1.75cm]{geometry}
\usepackage{amsmath}
\usepackage{graphicx}
\usepackage[colorlinks=true, allcolors=blue]{hyperref}
\usepackage[T1]{fontenc}
\usepackage[utf8]{inputenc}
\usepackage{amsmath, amssymb, amsthm}
\newenvironment{Proof*}[1][Proof]{%
  \par\noindent\textbf{#1.} \normalfont%
  }{\qed\par}

\newtheorem{definition}{Definition}
\newtheorem*{remark*}{Remark}

\newtheorem*{Example}{Example}
\usepackage{amsmath, amssymb}
\usepackage{tikz-cd}
\usepackage[utf8]{inputenc}
\usepackage{tikz}
\newtheorem{Proposition}{Proposition}
\newtheorem{theorem}{Theorem}
\newtheorem{corollary}{Corollary}
\usepackage{mathrsfs} 
\usepackage{hyperref}
\usepackage{geometry}
\usepackage{enumitem}
\theoremstyle{definition}
\title{Uniqueness of the non-commutative divergence cocycle}
\author{Pauline Baudat \footnote{ Department of Mathematics, University of Geneva, rue du Conseil Général 7-9, 1211 Geneva, Switzerland \textbf{e-mail : pauline.baudat@unige.ch}}}
\begin{document}
\maketitle
\tableofcontents
    \begin{abstract}
    We show that, for $n \ge 3 $, 1-cocycles of degree zero on the Lie algebra of derivations of the free associative algebra $T(A_n)$ with values in \( |T(A_n)| \otimes |T(A_n)| \) are linear combinations of the non-commutative divergence and its switch, when restricted to finite-degree quotients. Here, \( \rvert T(A_n) \rvert\) denotes the space of cyclic words. Furthermore, we study 1-cocycles of degree zero on the Lie algebra of symplectic derivations of the free Lie algebra $ \mathfrak{L}_{2n}$, and prove the uniqueness of the Enomoto–Satoh trace.
    \end{abstract}

\section{Introduction} 

Let 
$\operatorname{Poly}(\mathbb{R}^n) = \mathbb{R}[x_1, \dots, x_n]$
denote the space of polynomial functions over $\mathbb{R}$. Consider
\[
\mathfrak{X}(\mathbb{R}^n)
= \left\{
v = \sum_{i=1}^n v_i(x_1, \dots, x_n)\, \frac{\partial}{\partial x_i}
\ \middle|\
v_i \in \operatorname{Poly}(\mathbb{R}^n)
\right\},
\]
the \textbf{ Lie algebra of polynomial vector fields}, with the Lie bracket given by 
\[
[v,w] = v \circ w - w \circ v,
\qquad \forall\, v,w \in \mathfrak{X}(\mathbb{R}^n).
\]

 $\mathfrak{X} (\mathbb{R}^n)$ is graded by degree of the vector field. Indeed, each variable $x_i, \forall i \in \mathbb{N}$, is assigned a weight of 1. Thus, a homogeneous polynomial $P$ of degree 
$m$ is a linear combination of monomials of degree $m$. 

We call a map that sends every element of degree \(m\)  to an element of the same degree \(m\) a \textbf{degree zero} map.

A well known result is that the \textbf{ordinary divergence}
\begin{equation*}
\begin{aligned}
    \operatorname{div} : &\mathfrak{X} (\mathbb{R}^n) \to \operatorname{Poly}(\mathbb{R}^n) \\ 
    &v \mapsto \sum_{i=1}^n \frac{\partial_i}{\partial x_i}v_i ,
\end{aligned}
\end{equation*}
is the unique (up to scalar multiplication) degree zero function satisfying
 \begin{equation}\label{eq:divcocycle} 
         \operatorname{div} ([v_1,v_2]) = v_1( \operatorname{div} (v_2)) - v_2( \operatorname{div} (v_1)) \quad \forall v_1, v_2 \in \mathfrak{X} (\mathbb{R}^n). 
    \end{equation}

Functions satisfying \eqref{eq:divcocycle} are called \textbf{$1-$cocycles} and can be studied in more general setting. 
\\
In this paper, we will study the uniqueness of degree zero $1-$cocycles in the non-commutative case. Let $A_n$= $ \operatorname{Vect}_\mathbb{K}(x_1,...,x_n)$ the \(\mathbb{K}-\)vector space generated by $\{x_1,...,x_n\}$ with $\mathbb{K}$ a field of characteristic zero. The \textbf{free associative algebra} is \[T(A_n) =\bigoplus_{i \geq 1}^{\infty} A_n^{\otimes i}\]  the tensor algebra of $A_n$, which can be seen as the Lie algebra of the non-commutative polynomials with the bracket map 
\begin{equation*}
\begin{aligned} \relax 
    [X,Y] = X \otimes Y - Y \otimes X
&& \forall X,Y \in T(A_n).
\end{aligned}
\end{equation*}

A \textbf{derivation} of $T(A_n)$ is an endomorphism \(D\) of $T(A_n)$ satisfying 
 \begin{equation}\label{eq:1}
    D(X \otimes Y)= D(X) \otimes Y + X \otimes D(Y) \quad \forall X,Y \in T(A_n).
\end{equation}

In this setting, the set of all derivations of $T(A_n)$, denoted by $\operatorname{Der}(T(A_n))$, has a Lie algebra structure with the following bracket map  
\begin{equation}\label{eq:5}
[D_1,D_2]=D_1\circ D_2 - D_2 \circ D_1 \quad \forall D_1,D_2\in \operatorname{Der}(T(A_n)).
\end{equation}
There exists a direct sum decomposition
\begin{equation*}
    \operatorname{Der}(T(A_n)) \cong \bigoplus_{k \geq -1} \operatorname{Der}(T(A_n))(k),
\end{equation*}
where $\operatorname{Der}(T(A_n))(k)$ is the set of all derivation of degree k.
We also define $\operatorname{Der}_{\geq k }(T(A_n))$ the set of all derivations of degree larger than or equal to k. 

The quotient vector space
  \[
    \vert T(A_n)\vert := T(A_n)/ [T(A_n), 
    T(A_n)]
    \]
is called the \textbf{abelianization} of  $\operatorname{Der}(T(A_n))$. It is an abelian Lie algebra and can be seen as the space of cyclic words.

In our framework, we naturally introduce \textbf{the non-commutative divergence}, a linear map 
    \begin{equation} \label{Div}
        \operatorname{Div} : \operatorname{Der}(T(A_n)) \to |T(A_n)| \otimes |T(A_n)|, \quad \operatorname{Div}(u) = \sum_{i=1}^{n} |\partial_i u(z_i)|,
    \end{equation}
    where the \(i\)-th partial derivative is defined as 
    \begin{align*}
        \partial_i : \quad &T(A_n) \longrightarrow T(A_n) \otimes T(A_n) \\
                           &x_j \longmapsto \delta_{i,j} \, 1 \otimes 1,
    \end{align*}
and its \textbf{switch } 
    \(\sigma \circ \operatorname{Div} : \operatorname{Der}(T(A_n)) \to |T(A_n)| \otimes |T(A_n)|\), where \(\sigma\) is the transposition of the two elements of \(|T(A_n)|\). 

In \cite{Alekseev}, authors conjectured that every degree zero 1-cocycle on $\operatorname{Der}(T(A_n))$ with values in $\rvert T(A_n) \rvert \otimes \rvert T(A_n) \rvert$ is a linear combination of the maps Div and $\sigma \circ$ Div.
In an attempt to prove this conjecture, we use theorem \ref{thm:MoritaSakasai} from Shigeyuki Morita, Takuya Sakasai and Masaaki Suzuki \cite{Morita}. This gives us conditions on the generators of $\operatorname{Der}(T(A_n))$, which allows us to show our main result.
\begin{theorem}  \label{thm:mainresult}
    For all $n \ge 3$, let
    \begin{equation*}
    c : \operatorname{Der}(T(A_n)) \to \lvert T(A_n) \rvert \otimes \lvert T (A_n) \rvert \end{equation*} be a 1-cocycle of degree zero. Then, \(\forall k \leq n\), its restriction to the quotient by $\operatorname{Der}_{\geq k }(T(A_n))$ $\quad $
    \begin{equation*}
    c : \operatorname{Der}_{\geq 0}(T(A_n)) / \operatorname{Der}_{\geq k }(T(A_n)) \to \lvert T(A_n) \rvert \otimes \lvert T(A_n) \rvert, 
    \end{equation*}
    is a linear combination of the non-commutative divergence $\operatorname{Div}$ and $ \sigma \circ \operatorname{Div}$.  
\end{theorem}

This theorem gives us the following corollary. 

\begin{corollary}
    For all $ n \ge 3$, there exists no $c: \operatorname{Der}(T(A_n)) \to \rvert T(A_n) \rvert$ 1-cocycle of degree zero. 
\end{corollary}

Two remarks can be made about this statement: 
\begin{enumerate}
    \item The appearance of $ \operatorname{Der}_{\geq 0}(T(A_n)) / \operatorname{Der}_{\geq k }(T(A_n))$ is due to the condition in Morita-Sakasai-Suzuki's theorem. 
    \item The case of $n=1$ was studied by Alexander Gonsales \cite{n=1}. He proved that the space of 1-cocycles of degree 0 is 3-dimensional with the basis given by : $\operatorname{div} \otimes 1, 1 \otimes \operatorname{div}$ and $\operatorname{Div}$. An alternative proof of this fact will be given in section \ref{chapitre4}.
    \item The above result does not easily generalize to the case $n=2$. It is due to the fact that, in $T(A_2)$, we have $\rvert x_i x_j x_k \rvert =\rvert x_i x_k x_j \rvert \quad \forall i,k,j \in \{1,2\} $. More details can be found in remark \ref{n=2}. 
\end{enumerate}

The second part of this work deals with the specific case of the symplectic derivation Lie algebra of the free Lie algebra : \( \mathrm{Der}_{Sp}(\mathfrak{L}_{2n}) \). 

Let $\mathfrak{L}_{2n}$ \textbf{the free Lie algebra} generated by $\{x_1,...,x_n,y_1,...,y_n\}$ and $H=\operatorname{Vect}_\mathbb{K}(x_1,...,x_n,y_1,...,y_n)$.
A \textbf{derivation} of $\mathfrak{L}_{2n}$ is an endomorphisme D of $\mathfrak{L}_{2n}$ such that  
\[
D([z_{i_1},z_{i_2}])= [D(z_{i_1}),z_{i_2}] + [z_{i_1},D(z_{i_2})] \quad \forall z_{i_1},z_{i_2} \in \mathfrak{L}_{2n}
.\]
The \textbf{symplectic derivation Lie algebra of the free associative algebra} is 
    \begin{equation*}
        \operatorname{Der}_{Sp}(\mathfrak{L}_{2n})=\left\{ D \in \operatorname{Der}(\mathfrak{L}_{2n}) ; D\left(\sum_{j=1}^n [x_j,y_j]\right)=0\right\}.
    \end{equation*}
In this setting, we consider again $1$-cocylces of degree zero \( c: \mathrm{Der}_{Sp}(\mathfrak{L}_{2n}) \to |T(H)| \), up to scalar multiplication. A particular $1$-cocycle of degree zero was found in \cite{ES1} (see \eqref{def:ES}), and it is referred to as the Enomoto-Sato trace.  Our contribution is to show its uniqueness.

\begin{theorem} \label{thm:UniqueES}
     There is a unique, up to multiple, 1-cocycle of degree zero on $\operatorname{Der}_{Sp}(\mathfrak{L_{2n}})$ with values in $\rvert T(H) \rvert$. 
\end{theorem}

\section*{Acknowledgments}

I would like to express my sincere gratitude to Anton Alekseev for his guidance and support throughout this paper. I am also grateful to Vladimir Dotsenko for suggesting the article by Shigeyuki Morita, Takuya Sakasai and Masaaki Suzuki, as well as Takuya Sakasai, Takao Satoh and Nariya Kawazumi for taking the time to read this work.

\section{Preliminaries}

We will now recall some definitions of Lie algebras, derivations and 1-cocycles. This introduction follows the paper of Shigeyuki Morita, Takuya Sakasai, Masaaki Suzuki \cite{Morita}. 
\subsection{Lie algebra}

\begin{definition}
    A $\mathbb{K}$-vector space $\mathfrak{g}$  is called a \textbf{Lie algebra} if it has a $\mathbb{K}$-bilinear map 
\[
[ \cdot , \cdot ] : \mathfrak{g} \otimes \mathfrak{g} \to \mathfrak{g},
\]
satisfying the following properties 
\[
\begin{aligned}
    & 1.[X, Y] = -[Y, X] \quad && \forall X, Y \in \mathfrak{g} \quad && \text{(Antisymmetry)} \\
    & 2.[X, [Y, Z]] + [Y, [Z, X]] + [Z, [X, Y]] = 0 && \forall X, Y, Z \in \mathfrak{g}.\quad  && \text{(Jacobi identity)}
\end{aligned}
\]
The product $[X,Y]$ is called a Lie bracket.
\end{definition}

As defined above,  $T(A_n) =\bigoplus_{i \geq 1}^{\infty} A_n^{\otimes i}$ is the tensor algebra of $A_n$. It's the free associative Lie algebra with the following bracket map  
\begin{equation*}
\begin{aligned} \ 
    [X,Y] = X \otimes Y - Y \otimes X
&& \forall X,Y \in T(A_n)
\end{aligned}
\end{equation*}

$T(A_n)$ has a natural graduated structure given by the degree of the monomials. There exists a direct sum decomposition 
\begin{equation*}\label{eq:3}
T(A_n) = \bigoplus_{k \geq 1}^{\infty} T(A_n)(k)
\end{equation*}
with $T(A_n)(k)$ the vector space of homogeneous polynomials of degree $k$ and such that $[T(A_n)(k),T(A_n)(p)] \subset T(A_n)(k+p)$ $ \forall k,p \le 0.$

\begin{Example}
    \begin{equation*}
    \begin{aligned}
        &x_1 \otimes x_2 \in T(A_n)(2) \\
        &x_3 \otimes x_4 \otimes x_5 \in T(A_n)(3) \\
        &[x_1 \otimes x_2,x_3 \otimes x_4 \otimes x_5] = x_1 \otimes x_2 \otimes x_3 \otimes x_4 \otimes x_5 - x_3 \otimes x_4 \otimes x_5 \otimes x_1 \otimes x_2 \in T(A_n)(5)
    \end{aligned}
    \end{equation*}
\end{Example}
    
$\rvert T(A_n) \rvert$ is the abelianization of $T(A_n)$.
Indeed with $\pi$ the canonical projection on the abelianization, we have that the following diagram commutes
\[
\begin{tikzcd}
T(A_n) \otimes T(A_n) \arrow[r, "{[\cdot,\cdot]}"] \arrow[d, "\pi \otimes \pi"'] & T(A_n) \arrow[d, "\pi"] \\
{|T(A_n)| \otimes |T(A_n)|} \arrow[r, "0"] & {|T(A_n)|}.
\end{tikzcd}
\]
$\rvert T(A_n) \rvert$  can be viewed as the space of cyclic words.  

\begin{Example}
    \(\rvert x_1 x_2 x_3 \rvert = \rvert x_2 x_3 x_1 \rvert = \rvert x_3 x_1 x_2 \rvert \in  \rvert T(A_n) \rvert \)
    \text{but } \(\rvert x_1 x_2 x_3 \rvert \neq \rvert x_1 x_3 x_2 \rvert\).
\end{Example}

\subsection{Derivation of a Lie algebra}

$\operatorname{Der}(T(A_n))$, the set of all derivations of $T(A_n)$, is also graded by degree of derivation :  $D\in \operatorname{Der}(T(A_n))$ is a homogeneous element of degree \(d\) if it sends any  of degree $m$ element to a degree $d+m$ element. There exists a direct sum decomposition 
\begin{equation*}
    \operatorname{Der}(T(A_n)) \cong \bigoplus_{k \geq -1} \operatorname{Der}(T(A_n))(k),
\end{equation*}
where $\operatorname{Der}(T(A_n))(k)$ is the set of all derivation of degree $k$. 
\\
Let $A_n^*$ be the dual of $A_n$ with basis \{$x_1^*,x_2^*,...,x_n^*$\}. A derivation is uniquely determined by its action on the degree 1 component $T(A_n)(1) \cong A_n$. Conversely, every homomorphism in $\operatorname{Hom}(A_n,T(A_n))$ induces a derivation of $T(A_n)$. Hence, we obtain 
\begin{equation*}
\operatorname{Der}(T(A_n)) \cong \operatorname{Hom}(A_n, T(A_n)) 
\cong \bigoplus_{k \geq -1} \operatorname{Der}(T(A_n))(k),
\end{equation*}
where 
\begin{equation*}
\operatorname{Der}(T(A_n))(k) := \operatorname{Hom}(A_n, A_n^{\otimes (k+1)}) = A_n^* \otimes A_n^{\otimes (k+1)}
\end{equation*}
is the degree $k$ homogeneous part of $\operatorname{Hom}(A_n,T(A_n))$. As for $T(A_n)$, we have $[\operatorname{Der}(T(A_n))(p), \operatorname{Der}(T(A_n))(q)] \subset \operatorname{Der}(T(A_n))(p+q), \,\, \forall p,q \geq -1$.
\begin{Example}
    We can express
    \begin{equation*}
        \begin{aligned}
            D: A_2 &\to T(A_2) \\
            x_1 &\mapsto x_1 \otimes x_2 \\
            x_2 &\mapsto x_2 \otimes x_1
        \end{aligned}
    \end{equation*}
    as
    \begin{equation*}
        D= x_1^* \otimes x_1 \otimes x_2 + x^*_2 \otimes x_2 \otimes x_1 \in \operatorname{Der}(T(A_n))(1) \cong A_n^* \otimes A_n^2.
    \end{equation*}
\end{Example}

\begin{remark*}
    For two elements 
\begin{align*}
D_1 &= f \otimes u_1 \otimes u_2 \otimes \cdots \otimes u_{p+1} \in  A_n^* \otimes A_n^{\otimes (p+1)},\\
D_2 &= g \otimes v_1 \otimes v_2 \otimes \cdots \otimes v_{q+1} \in A_n^* \otimes A_n^{\otimes (q+1)},
\end{align*}
where \( f, g \in A_n^* \) and \( u_1, \ldots, u_{p+1}, v_1, \ldots, v_{q+1} \in A_n \), the bracket map \eqref{eq:5} \([D_1, D_2] \in  A_n^*\otimes A_n^{\otimes (p+q+1)}\) can be defined as
\begin{equation*}
\begin{aligned} \
[D_1, D_2] =
& \sum_{s=1}^{q+1} f(v_s)g \otimes v_1 \otimes \cdots \otimes v_{s-1} \otimes (u_1 \otimes \cdots \otimes u_{p+1}) \otimes v_{s+1} \otimes \cdots \otimes v_{q+1} \\
& - \sum_{t=1}^{p+1} g(u_t)f \otimes u_1 \otimes \cdots \otimes u_{t-1} \otimes (v_1 \otimes \cdots \otimes v_{q+1}) \otimes u_{t+1} \otimes \cdots \otimes u_{p+1}.
\end{aligned}
\end{equation*}
\end{remark*}

\begin{remark*}
    Note that $\operatorname{Der}_{\geq 0}(T(A_n)) / \operatorname{Der}_{\geq k }(T(A_n)) $ represents the set of all derivation of non-negative degree smaller or equal to k. 
\end{remark*}

\subsection{1-cocycles and non-commutative calculus}

We define a 1-cocycle as follows  

\begin{definition}
    Let $\mathfrak{g}$ be a Lie algebra and $M$ a $\mathfrak{g}$-module.  
    A linear map $c : \mathfrak{g} \to M$ is a \textbf{1-cocycle of Lie algebra} with values in M if it satisfies  
    \begin{equation} \label{cocycle}
        \forall u, v \in \mathfrak{g}, \quad c([u, v]) = u(c(v)) - v(c(u)).
    \end{equation}
\end{definition}

\begin{remark*}
    This definition follows from the cocycle condition $dc=0$, with  
    \begin{equation} \label{eq:dc}
    (dc)(x_1, \dots, x_n) = \sum_{1 \leq i < j \leq n} (-1)^{i+j+1} c([x_i, x_j], x_1, \dots, \widehat{x_i}, \dots, \widehat{x_j}, \dots, x_n) + 
    \sum_{i=1}^{n} (-1)^i x_i c(x_1, \dots, \widehat{x_i}, \dots, x_n).
    \end{equation}
    for all $ x_1,\dots,x_n \in \operatorname{Der}(T(A_n))$
\end{remark*}


\begin{Example}
We can consider the following example for the maps \(\operatorname{Div}\) and $\sigma \circ\operatorname{Div} $ defined in the introduction \eqref{Div}
    \begin{align*}
        &\operatorname{Div} (x_{i_0}^* \otimes x_{i_1} \otimes \cdots \otimes x_{i_{k+1}}) = \sum_{j=1}^{k+1} \delta_{i_o,i_j} \rvert x_{i_1} \cdots x_{i_{j-1}} \rvert \otimes \rvert x_{j_{j+1}} \cdots x_{i_k} \rvert, \\
        &\sigma \circ \operatorname{Div}(x_{i_0}^* \otimes x_{i_1} \otimes \cdots  \otimes x_{i_{k+1}}) = \sum_{j=1}^{k+1} \delta_{i_o,i_j} \rvert x_{j_{j+1}} \cdots x_{i_k} \rvert \otimes\rvert x_{i_1} \cdots x_{i_{j-1}} \rvert     . 
    \end{align*}
    
\end{Example}

\section{Main result}

\subsection{The generators of \(\operatorname{Der}(T(A_n))\)}

As mentioned in the introduction, to prove the main theorem of this paper, we will use the following result given by Shigeyuki Morita, Takuya Sakasai and Masaaki Suzuki \cite{Morita}.

\begin{theorem}\label{thm:MoritaSakasai}
    For n $\geq$ 2, we have a direct sum decomposition 
    \begin{equation*}
    \operatorname{Der}(T(A_n))(2) = A_n^{\otimes 2} \oplus [\operatorname{Der}(T(A_n))(1), \operatorname{Der}(T(A_n))(1)].
    \end{equation*}
    
    If $n \geq k \geq 3$ , we have
    \begin{equation*}
    \operatorname{Der}(T(A_n))(k)
    = [\operatorname{Der}(T(A_n))(k-1), \operatorname{Der}(T(A_n))(1)]
    + [\operatorname{Der}(T(A_n))(k-2), \operatorname{Der}(T(A_n))(2)].
    \end{equation*}
\end{theorem}

\begin{remark*}
    This theorem demonstrates that $\forall n\geq 3,  \operatorname{Der}(T(A_n))(k)$ where $k \leq n$ is completely determined by $\operatorname{Der}(T(A_n))(1)$ and $\operatorname{Der}(T(A_n))(2)$. This observation will be crucial for the proof of theorem \ref{thm:mainresult}.
\end{remark*}

\begin{remark*}
    In this theorem, $A_n^{\otimes 2}$ is seen as domain of the injective function \( s: A_n^{\otimes 2} \to \operatorname{Der}(T(A_n))(2) \) such that \(s(x_i \otimes x_j)= x_1 ^* \otimes x_i \otimes x_1 \otimes x_j\). Indeed, we have that $s(A_n^{\otimes 2})$ and elements of the form $[\operatorname{Der}(T(A_n))(1), \operatorname{Der}(T(A_n))(1)]$ generate $\operatorname{Der}(T(A_n))(2)$.
\end{remark*}

\subsection{Proof of the main result}

In this section, we will prove theorem \ref{thm:mainresult}. The proof is divided into 4 parts :
\begin{enumerate}[label=\Roman*.]
    \item The restriction of c to \( \operatorname{Der}(T(A_n))(0)\)
    \item The restriction of c to \( \operatorname{Der}(T(A_n))(1)\)
    \item The restriction of c to \( \operatorname{Der}(T(A_n))(2)\)
    \item The restriction of c to \( \operatorname{Der}(T(A_n))(3)\) and the conclusion.
\end{enumerate}

\begin{remark*}
For the proof, we will use the fact that $c$ is $\operatorname{Gl}_n\mathbb{K}$- equivariant. This result is based on the following observations : 
\begin{enumerate}
    \item Authors \cite{Morita} showed that $\operatorname{Der}(T(A_n))(0) \cong \mathfrak{gl}_n\mathbb{K}$.
    \item $c(A)= a \operatorname{Tr}(A)$ for  $A\in \mathfrak{gl}_n\mathbb{K}$ (see \eqref{tr}).  Since $c$ is a cocycle, we have 
     \begin{equation*}
     \begin{aligned}
         c([A,D])&=A \cdot c(D) - D \cdot c(A) \\
         &=  A \cdot c(D) - D \cdot a\operatorname{Tr}(A) \\
         &= A \cdot c(D) 
         \quad \forall A \in \mathfrak{gl}_n\mathbb{K} \quad  \forall D \in \operatorname{Der}(T(A_n)).
     \end{aligned}
     \end{equation*}
     In other words, $c$ is $\mathfrak{gl}_n\mathbb{K}$-equivariant. 
     \item  We have these group actions  
     \begin{enumerate}
        \item $\operatorname{Gl}_n\mathbb{K} \curvearrowright A_n $ 
        \begin{equation*}
            X \cdot x_{i_1} = \sum _{l=1}^n X_{i_1l}x_{l} \quad \forall X\in  \operatorname{Gl}_n\mathbb{K}, \quad \forall x_{i_1} \in A_n.
        \end{equation*}
        \item $\operatorname{Gl}_n\mathbb{K} \curvearrowright A_n^* $  
        \begin{equation*}
            X \cdot  x_{i_0}^*(x_j) = x_{i_0}^*(X^{-1}\cdot x_{i_1}) \quad \forall X\in  \operatorname{Gl}_n\mathbb{K}, \quad \forall x_{i_0}^* \in A_n^*, \quad \forall x_{i_1} \in A_n.
        \end{equation*}
        \item $\operatorname{Gl}_n\mathbb{K} \curvearrowright \rvert T(A_n) \rvert \otimes \rvert T(A_n) \rvert $  
        \begin{equation*}
        \begin{aligned}
            &X \cdot  (\rvert x_{i_1} \cdots x_{i_l} \rvert \otimes \rvert x_{i_{l+1}}\cdots x_{i_k} \rvert  ) = \rvert X \cdot x_{i_1}\cdots X \cdot x_{i_l} \rvert \otimes \rvert X \cdot x_{i_{l+1}} \cdots X \cdot x_{i_k} \rvert  \\
            &\forall X\in  \operatorname{Gl}_n\mathbb{K}, \quad \forall \rvert x_{i_1} \cdots x_{i_l} \rvert, \rvert x_{i_{l+1}} \cdots x_{i_k} \rvert \in \rvert T(A_n) \rvert .
        \end{aligned}
        \end{equation*}
        \item $\operatorname{Gl}_n\mathbb{K} \curvearrowright \operatorname{Der}(T(A_n)) $
        \begin{equation*}
        X \cdot (x_{i_0}^* \otimes x_{i_1} \otimes \cdots \otimes x_{i_{k+1}}) = X \cdot x_{i_0}^* \otimes X \cdot x_{i_1} \otimes \cdots \otimes X \cdot x_{i_{k+1}} \quad \forall X\in  Gl_n\mathbb{K} , \quad \forall x_{i_0}^* \otimes x_{i_1} \otimes \cdots \otimes x_{i_{k+1}} \in \operatorname{Der}(T(A_n)).
        \end{equation*}
    \end{enumerate}
    
\end{enumerate}
    Combining those facts, $c$ is $\operatorname{Gl}_n\mathbb{K}$-equivariant by canonical integration.
    In particular, it means that $c$ can always be written as 
\begin{equation} \label{equiv}
    c(x_{i_0}^* \otimes x_{i_1} \otimes \cdots \otimes x_{i_{k+1}}) = \sum_{j=1}^{k+1} c_j x_{i_0} ^*(x_{i_j}) \rvert x_{i_1} \cdots x_{i_{j-1}} \rvert \otimes \rvert x_{i_{j+1}} \cdots x_{i_{k+1}} \rvert.
\end{equation}

\end{remark*}

\begin{Proof*} 
I. As mentioned earlier $\operatorname{Der}(T(A_n))(0) \cong \mathfrak {gl}_n\mathbb{K}$. Since c is a cocycle, with \eqref{eq:dc}, we have  
\begin{equation*}
dc(A,B)=-c[A,B]=0 \quad \forall A,B \in \mathfrak {gl}_n\mathbb{K}.
\end{equation*}
Which means that $c(AB)=c(BA)$. 
Therefore, we have, for \(E_{ij}=\mathbf{1}_{\{i=j\}}\),
\begin{equation} \label{tr}
\begin{aligned}
    &c(E_{ij}E_{jk})=c(E_{ik})=c(E_{jk}E_{ij})=0 \quad \forall i\ne j,   \\
    &c(E_{ii})=c(E_{ik}E_{ki})=c(E_{ki}E_{ik})=c(E_{kk}) \quad \forall 1\leq i,k \leq n.
\end{aligned}
\end{equation}
This implies that $c(E_{ii})=c(E_{jj}) = \frac{c(Id)}{n} \; \; \forall 1 \le i \le n$.
For the rest of the proof, we set $a=\frac{c(Id)}{n} \in \mathbb{K} $, which means that we have 

\begin{equation*}
    c(x_{i_0}^* \otimes x_{i_1}) = \delta_{{i_o},{i_1}} a.
\end{equation*}

\vspace{1em} 

II.  With the remark \eqref{equiv}, we know $\forall x_{j_0},x_{j_1},x_{j_2} \in A_n, \; \exists \alpha, \beta, \gamma,\omega \in \mathbb{K}$ such that 

\begin{equation*}
c(x^*_{j_0} \otimes x_{j_1} \otimes x_{j_2}) = 
\delta_{j_0,j_1} \big( 
\alpha \lvert x_{j_2} \rvert \otimes 1 
+ \beta \, 1 \otimes \lvert x_{j_2} \rvert 
\big) 
+ \delta_{j_0,j_2} \big( 
\gamma \lvert x_{j_1} \rvert \otimes 1 
+ \omega \, 1 \otimes \lvert x_{j_1} \rvert 
\big).
\end{equation*}
To determine them, we are using the following relation  

\begin{equation*}
[x^*_{i_0}, x^*_{j_0} \otimes x_{j_1} \otimes x_{j_2}] 
= x^*_{i_0}(x_{j_1}) x^*_{j_0} \otimes x_{j_2} 
+ x^*_{i_0}(x_{j_2}) x^*_{j_0} \otimes x_{j_1}.
\end{equation*}
We find 

\begin{equation*}
\begin{aligned}
    LHS &= c([x^*_{i_0}, x^*_{j_0} \otimes x_{j_1} \otimes x_{j_2}]) \\
    &= x_{i_0}^*(c(x^*_{j_0} \otimes x_{j_1} \otimes x_{j_2})) - x^*_{j_0} \otimes x_{j_1} \otimes x_{j_2}(c(x^*_{i_0}))\\
    &= x_{i_0}^*\Big( 
    \delta_{j_0,j_1} \big( 
    \alpha \lvert x_{j_2} \rvert \otimes 1 + \beta \, 1 \otimes \lvert x_{j_2} \rvert \big) 
    + \delta_{j_0,j_2} \big( 
    \gamma \lvert x_{j_1} \rvert \otimes 1 + \omega \, 1 \otimes \lvert x_{j_1} \rvert \big) 
    \Big) - 0  \\
    &= \delta_{j_0,j_1}  \delta_{i_0,j_2} (\alpha + \beta) 1 \otimes 1  
    + \delta_{j_0,j_2}  \delta_{i_0,j_1} (\gamma + \omega) 1 \otimes 1, \\
    \\
    RHS &= c([x^*_{i_0}, x^*_{j_0} \otimes x_{j_1} \otimes x_{j_2}]) \\ 
    &= c(x^*_{i_0}(x_{j_1}) x^*_{j_0} \otimes x_{j_2} 
    + x^*_{i_0}(x_{j_2}) x^*_{j_0} \otimes x_{j_1}) \\
    &= \delta_{i_0,j_1} c( x^*_{j_0} \otimes x_{j_2}) 
    + \delta_{i_0,j_2} c(x^*_{j_0} \otimes x_{j_1}) \\
    &= \delta_{i_0,j_1} \delta_{j_0,j_2} a 1 \otimes 1 
    + \delta_{i_0,j_2} \delta_{j_0,j_1} a 1 \otimes 1.
\end{aligned}
\end{equation*}
Combining these equalities we have

\begin{equation} \label{eq: degree 1}
\left\{
\begin{aligned}
    \alpha + \beta &= a \\
    \gamma + \omega &= a.
\end{aligned}
\right.
\end{equation}

\vspace{1em} 

III. As before, with \eqref{equiv}, we know $ \exists a_1,a_2,a_3,a_4,b_1,b_2,b_3,b_4,c_1,c_2,c_3,c_4 \in \mathbb{K}$ such that 

\begin{equation*}
\begin{aligned}
    c(x^*_{j_0} \otimes x_{j_1} \otimes x_{j_2} \otimes x_{j_3}) = 
    &\ \delta_{j_0,j_1} \big( 
    a_1 \lvert x_{j_2} x_{j_3} \rvert \otimes 1 
    + a_2 \, 1 \otimes \lvert x_{j_2} x_{j_3} \rvert 
    + a_3 \lvert x_{j_2} \rvert \otimes \lvert x_{j_3} \rvert 
    + a_4 \lvert x_{j_3} \rvert \otimes \lvert x_{j_2} \rvert 
    \big) \\
    &+ \delta_{j_0,j_2} \big( 
    b_1 \lvert x_{j_1} x_{j_3} \rvert \otimes 1 
    + b_2 \, 1 \otimes \lvert x_{j_1} x_{j_3} \rvert 
    + b_3 \lvert x_{j_1} \rvert \otimes \lvert x_{j_3} \rvert 
    + b_4 \lvert x_{j_3} \rvert \otimes \lvert x_{j_1} \rvert 
    \big) \\
    &+ \delta_{j_0,j_3} \big( 
    c_1 \lvert x_{j_1} x_{j_2} \rvert \otimes 1 
    + c_2 \, 1 \otimes \lvert x_{j_1} x_{j_2} \rvert 
    + c_3 \lvert x_{j_1} \rvert \otimes \lvert x_{j_2} \rvert 
    + c_4 \lvert x_{j_2} \rvert \otimes \lvert x_{j_1} \rvert
    \big).
\end{aligned}
\end{equation*}
In order to identify them, we consider the following equalities 

\begin{align*}
     \begin{aligned}
        1. & x^*_{l} \otimes x_{l} \otimes x_{i_1} \otimes x_{i_2}=[x^*_{i_1} \otimes x_{i_1} \otimes x_{i_2}, x^*_{l} \otimes x_{l} \otimes x_{i_1}], \\  
        2. &x^*_{l} \otimes x_{i_1} \otimes x_{i_2} \otimes x_{l}=[x^*_{i_1} \otimes x_{i_1} \otimes x_{i_2}, x^*_{l} \otimes x_{i_1} \otimes x_{l}], \\  
        3. &[x^*_{i_0}, x^*_{l} \otimes x_{i_1} \otimes x_{l} \otimes x_{i_2}] 
        = \delta_{i_0,i_1} x^*_{l} \otimes x_{l} \otimes x_{i_2} 
        + \delta_{i_0,l} x^*_{l} \otimes x_{i_1} \otimes x_{i_2} 
        + \delta_{i_0,i_2} x^*_{l} \otimes x_{i_1} \otimes x_{l}.
    \end{aligned}
\end{align*}

1. For the first equality we have 
\begin{align*}
    \text{LHS} &= c(x^*_{l} \otimes x_{l} \otimes x_{i_1} \otimes x_{i_2}) \\ 
    &= a_1 \lvert x_{i_1} x_{i_2} \rvert \otimes 1  
    + a_2 \, 1 \otimes \lvert x_{i_1} x_{i_2} \rvert  
    + a_3 \lvert x_{i_1} \rvert \otimes \lvert x_{i_2} \rvert  
    + a_4 \lvert x_{i_2} \rvert \otimes \lvert x_{i_1} \rvert. \\
    \\
    \text{RHS} &= c\big([x^*_{i_1} \otimes x_{i_1} \otimes x_{i_2}, x^*_{l} \otimes x_{l} \otimes x_{i_1}]\big) \\ 
    &= x^*_{i_1} \otimes x_{i_1} \otimes x_{i_2} (c(x^*_{l} \otimes x_{l} \otimes x_{i_1}) )- x^*_{l} \otimes x_{l} \otimes x_{i_1}(c(x^*_{i_1} \otimes x_{i_1} \otimes x_{i_2})) \\ 
    &= x^*_{i_1} \otimes x_{i_1} \otimes x_{i_2} \cdot (\alpha \lvert x_{i_1} \rvert \otimes 1 + \beta 1 \otimes \lvert x_{i_1} \rvert) -0 \\ 
    &= \alpha \lvert x_{i_1} x_{i_2} \rvert \otimes 1 + \beta 1 \otimes \lvert x_{i_1} x_{i_2} \rvert.
\end{align*}
Which gives us  
 \begin{equation} \label{eq : a}
\left\{
\begin{aligned}
    &a_1 = \alpha \\
    &a_2 = \beta  \\
    &a_3 =0 \\
    &a_4=0.\\
\end{aligned}
\right.
\end{equation}

2. To find $c_1,c_2,c_3,c_4$ we have 
\begin{align*}
    \text{LHS} & = c(x^*_{l} \otimes x_{i_1} \otimes x_{i_2} \otimes x_{l}) \\
    &= \big(c_1 \lvert x_{i_1} x_{i_2} \rvert \otimes 1 
    + c_2 \, 1 \otimes \lvert x_{i_1} x_{i_2} \rvert 
    + c_3 \lvert x_{i_1} \rvert \otimes \lvert x_{i_2} \rvert 
    + c_4 \lvert x_{i_2} \rvert \otimes \lvert x_{i_1} \rvert 
    \big), \\
    \\
    \text{RHS} &= c([x^*_{i_1} \otimes x_{i_1} \otimes x_{i_2}, x^*_{l} \otimes x_{i_1} \otimes x_{l}]) \\ 
    &= x^*_{i_1} \otimes x_{i_1} \otimes x_{i_2}(c(x^*_{l} \otimes x_{i_1} \otimes x_{l})) - x^*_{l} \otimes x_{i_1} \otimes x_{l}(c(x^*_{i_1} \otimes x_{i_1} \otimes x_{i_2})) \\
    &=x^*_{i_1} \otimes x_{i_1} \otimes x_{i_2}( \gamma \rvert x_{i_1} \rvert \otimes 1 + \omega 1 \otimes \rvert x_{i_1} \rvert) - 0\\
    &=\gamma \rvert x_{i_1}  x_{i_2} \rvert \otimes 1 + \omega 1 \otimes \rvert x_{i_1}  x_{i_2} \rvert.
\end{align*}
As LHS equals RHS, we find 
\begin{equation}  \label{eq : c}
\left\{
\begin{aligned}
    &c_1 = \gamma \\
    &c_2 = \omega  \\
    &c_3 =0 \\
    &c_4=0.\\
\end{aligned}
\right.
\end{equation}

3. This time we have  

\begin{align*}
    RHS &= c( \delta_{i_0,i_1} x^*_{l} \otimes x_{l} \otimes x_{i_2} 
      + \delta_{i_0,i_l} x^*_{l} \otimes x_{i_1} \otimes x_{i_2} 
      + \delta_{i_0,i_2} x^*_{l} \otimes x_{i_1} \otimes x_{i_2}) \\ 
    &= \delta_{i_0,i_1} (\alpha \rvert x_{i_2} \rvert \otimes 1 
      + \beta 1 \otimes \rvert x_{i_2} \rvert) + \delta_{i_0,i_2} (\gamma \rvert x_{i_1} \rvert \otimes 1 
      + \omega 1 \otimes \rvert x_{i_1} \rvert), \\
      \\
    LHS &=c([x^*_{i_0}, x^*_{l} \otimes x_{i_1} \otimes x_{l} \otimes x_{i_2}]) \\
    &= x^*_{i_0}(c(x^*_{l} \otimes x_{i_1} \otimes x_{l} \otimes x_{i_2})) - x^*_{l} \otimes x_{i_1} \otimes x_{l} \otimes x_{i_2}(c(x^*_{i_0}))\\
    &= x^*_{i_0}\big(b_1 \lvert x_{i_1} x_{i_2} \rvert \otimes 1 
    + b_2 \, 1 \otimes \lvert x_{i_1} x_{i_2} \rvert
    + b_3 \lvert x_{i_1} \rvert \otimes \lvert x_{i_2} \rvert 
    + b_4 \lvert x_{i_2} \rvert \otimes \lvert x_{i_1} \rvert\big) -0 \\
    &= b_1 \big(\delta_{i_0,i_1} \rvert x_{i_2} \rvert \otimes 1 
    + \delta_{i_0,i_2} \rvert x_{i_1} \rvert \otimes 1\big) + b_2 \big(\delta_{i_0,i_1} 1 \otimes \rvert x_{i_2} \rvert  
    + \delta_{i_0,i_2} 1 \otimes \rvert x_{i_1} \rvert\big) \\
    &\quad + b_3 \big(\delta_{i_0,i_1} 1 \otimes \rvert x_{i_2} \rvert  
    + \delta_{i_0,i_2} \rvert x_{i_1} \rvert \otimes 1\big) + b_4 \big(\delta_{i_0,i_1} \rvert x_{i_2} \rvert \otimes 1 
    + \delta_{i_0,i_2} 1 \otimes \rvert x_{i_1} \rvert\big).
\end{align*}
Similarly, we have the following equalities $\forall t \in \mathbb{K} $ 
\begin{equation}  \label{eq : b}
\left\{
\begin{aligned}
    b_1 & = \alpha - t \\
    b_2 & = \beta - \gamma + \alpha - t \\
    b_3 & = \gamma - \alpha + t \\
    b_4 & = t.
\end{aligned}
\right.
\end{equation}

\vspace{1em} 

IV. With theorem \ref{thm:MoritaSakasai}, we find  

\begin{align*}
    c(x_{i_0}^* \otimes x_{i_1} \otimes x_{i_2} \otimes x_{i_3} \otimes x_{i_4})&= \delta_{i_0,i_1} (\alpha \rvert x_{i_2}x_{i_3}x_{i_4} \rvert \otimes 1) + \beta 1 \otimes \rvert x_{i_2}x_{i_3}x_{i_4} \rvert) \\
    &+\delta_{i_0,i_2} (b_1 \rvert x_{i_1}x_{i_3}x_{i_4} \rvert \otimes 1) + b_2 1 \otimes \rvert x_{i_1}x_{i_3}x_{i_4} \rvert+ b_3 \rvert x_{i_1} \rvert \otimes \rvert x_{i_3}x_{i_4}\rvert + b_4 \rvert x_{i_3}x_{i_4}\rvert \otimes \rvert x_{i_1} \rvert) \\
    &+\delta_{i_0,i_3} (b_1 \rvert x_{i_1}x_{i_2}x_{i_4} \rvert \otimes 1) + b_2 1 \otimes \rvert x_{i_1}x_{i_2}x_{i_4} \rvert+ b_3 \rvert x_{i_1}x_{i_2} \rvert \otimes \rvert x_{i_4}\rvert + b_4 \rvert x_{i_4}\rvert \otimes \rvert x_{i_1}x_{i_2} \rvert) \\
    &+ \delta_{i_0,i_4} (\gamma \rvert x_{i_1}x_{i_2}x_{i_3} \rvert \otimes 1) + \omega 1 \otimes \rvert x_{i_1}x_{i_2}x_{i_3} \rvert).
\end{align*}
Additional details can be found at the end of this paper in Appendix \ref{appendix}. 

We also know that 
\begin{align*}
    c([x_{i_0}^* \otimes x_{i_1} \otimes x_{i_2},x_{j_0}^* \otimes x_{j_1} \otimes x_{j_2} \otimes x_{j_3}]) =& \delta_{i_0,j_1} c(x_{j_0}^* \otimes x_{i_1} \otimes x_{i_2} \otimes x_{j_2} \otimes x_{j_3}) + \delta_{i_0,j_2} c(x_{j_0}^* \otimes x_{j_1} \otimes x_{i_1} \otimes x_{i_2} \otimes x_{j_3}) \\
    & + \delta_{i_0,j_3} c(x_{j_0}^* \otimes x_{j_1} \otimes x_{j_2} \otimes x_{i_1} \otimes x_{i_2}) - \delta_{j_0,i_1} c(x_{i_0}^* \otimes x_{j_1} \otimes x_{j_2} \otimes x_{j_3} \otimes x_{i_2}) \\
    &- \delta_{j_0,i_2} c(x_{i_0}^* \otimes x_{i_1} \otimes x_{j_1} \otimes x_{j_2} \otimes x_{j_3}).     
\end{align*}
Since c is a 1-cocycle we have in  $\delta_{i_0,j_2} \delta_{j_0,i_2}$

\begin{align*}
    LHS = & 0, \\
    RHS = & b_1 \rvert x_{j_1} x_{i_1} x_{j_3} \rvert \otimes 1 + b_2 1 \otimes \rvert x_{j_1} x_{i_1} x_{j_3} \rvert + b_3 \rvert x_{j_1} x_{i_1} \rvert \otimes \rvert x_{j_3} \rvert + b_4 \rvert x_{j_3} \rvert \otimes \rvert x_{j_1} x_{i_1} \rvert \\
    & - b_1 \rvert x_{i_1} x_{j_1} x_{j_3} \rvert \otimes 1 - b_2 1 \otimes \rvert x_{i_1} x_{j_1} x_{j_3} \rvert - b_3 \rvert x_{i_1} x_{j_1} \rvert \otimes \rvert x_{j_3} \rvert - b_4 \rvert x_{j_3} \rvert \otimes \rvert x_{i_1} x_{j_1} \rvert.
\end{align*}
Which gives us the last condition 
\begin{equation} \label{eq : b2}
 \quad b_1=b_2=0.
\end{equation}
By gathering \eqref{eq: degree 1},\eqref{eq : a}, \eqref{eq : c}, \eqref{eq : b}, \eqref{eq : b2}, we obtain  

\begin{align*}
    c(x^*_{i_0} \otimes x_{i_1} \otimes x_{i_2} \otimes x_{i_3} \otimes x_{i_4}) & = a( \delta_{i_0,i_1} (\alpha \rvert x_{i_2}x_{i_3}x_{i_4} \rvert \otimes 1 + (a-\alpha) 1 \otimes \rvert x_{i_2}x_{i_3}x_{i_4} \rvert) \\ &+ \delta_{i_0,i_2} ((a-\alpha) \rvert x_{i_1} \rvert \otimes \rvert x_{i_3}x_{i_4} \rvert + \alpha \rvert x_{i_3}x_{i_4} \rvert \otimes \rvert x_{i_1} \rvert) \\
    & + \delta_{i_0,i_3} ( (a-\alpha) \rvert x_{i_1}x_{i_2} \rvert \otimes \rvert x_{i_4} \rvert + \alpha \rvert x_{i_4} \rvert \otimes \rvert x_{i_1}x_{i_2} \rvert ) \\ & + \delta_{i_o,i_4} (\alpha \rvert x_{i_1}x_{i_2}x_{i_3} \rvert \otimes 1 + (a-\alpha) 1 \otimes \rvert x_{i_1}x_{i_2}x_{i_3} \rvert)) \\
    &= a\alpha \operatorname{Div}(x^*_{i_0} \otimes x_{i_1} \otimes x_{i_2} \otimes x_{i_3} \otimes x_{i_4}) + a(a- \alpha) \sigma \circ \operatorname{Div} (x^*_{i_0} \otimes x_{i_1} \otimes x_{i_2} \otimes x_{i_3} \otimes x_{i_4}).
\end{align*}
  To finish the proof, we need to state that, by Theorem \ref{thm:MoritaSakasai}, $c: \operatorname{Der}(T(A_n))(k) \to \rvert T(A_n) \rvert \otimes \rvert T(A_n) \rvert \; \; \forall k \ge 3 $ is completely determined by the lower degrees. 

\end{Proof*}

\begin{remark*} \label{n=2}
    The case \(n=2\) cannot be studied in the same way. The proof remains valid for parts I, II, and III, but the last step does not hold.
\end{remark*}

\begin{corollary}
    For all $ n \ge 3$, there exists no $c: \operatorname{Der}(T(A_n)) \to \rvert T(A_n) \rvert$ 1-cocycle of degree zero. 
\end{corollary}

\begin{Proof*}
    Suppose there $\exists \overline c: \operatorname{Der}(T(A_n)) \to \rvert T(A_n) \rvert$ a 1-cocycle of degree zero. 
    
    This implies that we can define $c_1(D) = 1 \otimes \overline c(D) \; \; \forall D \in \operatorname{Der}(T(A_n)) $ a 1-cocycle of degree zero.
    \\
    But there exist no $\alpha,a \in \mathbb{K}$ such that $c_1 : \operatorname{Der}(T(A_n)) \to \rvert T(A_n) \rvert \otimes \rvert T(A_n) \rvert$ is a 1-cocycle of degree zero.

 \end{Proof*}

\subsection{Case $n=1$} \label{chapitre4}

The case $n=1$ differs from \(n \ge 3\) because $\mathbb{K}\left<x\right>$ is commutative. We can consider $\{x^k|k \in \mathbb{Z}_{\geq 0}\}$ the basis of $\mathbb{K}\left<x\right>$, and  $\{x^k \otimes x^l|k,l \in \mathbb{Z}_{\geq 0}\}$ a basis of $\mathbb{K}\left<x\right> \otimes \mathbb{K}\left<x\right>$. 

As $\mathbb{K}\left<x\right>$ is commutative, we can consider $\{x^* \otimes x^k|k \in \mathbb{Z}_{\geq 0}\}$ as a basis of \(\operatorname{Der}(\mathbb{K}\left<x\right>)\). Note that $x^* \otimes x^{k+1}$ is an element of degree k in \(\operatorname{Der}(\mathbb{K}\left<x\right>)\). 

We want to know how many 1-cocycle $c : \operatorname{Der}(\mathbb{K}\left<x\right>) \to \mathbb{K}\left<x\right> \otimes \mathbb{K}\left<x\right>$ of degree zero exist. A first observation is  \[ \sigma \circ \operatorname{Div} = \operatorname{Div} \] due to the commutativity of $\mathbb{K}\left<x\right>$. Another observation is  \(div \otimes 1 \) and \(1 \otimes div\) are 1-cocycles of degree zero for $\mathbb{K}\left<x\right> \otimes  \mathbb{K}\left<x\right>$.

In his master thesis, Alexander Gonsales \cite{n=1} found the following result.


\begin{Proposition}
    The space of the 1-cocycles of degree zero from \(\operatorname{Der}(\mathbb{K}\left<x\right>)\) in $\mathbb{K}\left<x\right> \otimes \mathbb{K}\left<x\right>$ has dimension 3, i.e $\operatorname{dim}(H_{(0)}^1(\operatorname{Der}(\mathbb{K}\left<x\right>),\mathbb{K}\left<x\right> \otimes \mathbb{K}\left<x\right>))=3$. 
    We can introduce a basis such that the three cocycles are 
 \begin{enumerate}
        \item commutative left-divergence $\operatorname{div} \otimes 1$,
        \item commutative right-divergence $1 \otimes \operatorname{div}$,
        \item non-commutative divergence.
    \end{enumerate}
\end{Proposition}

\begin{remark*}
    For the proof of the proposition, it is important to note that
    \begin{equation*} 
    \begin{aligned} \
        [x^* \otimes x^{k+1}, x^* \otimes x^{l+1}]&= x^* \otimes x^{k+1}(x^* \otimes x^{l+1})- x^* \otimes x^{l+1}(x^* \otimes x^{k+1}) \\
        &=(l+1) x^* \otimes x^{k+l+1} - (k+1) x^* \otimes x^{k+l+1} \\
        &=(l-k)x^* \otimes x^{k+l+1}.
    \end{aligned}
    \end{equation*}
\end{remark*}

\begin{Proof*}
    Any 1-cocycles of degree zero c: $\operatorname{Der}(\mathbb{K}\left<x\right>) \to \mathbb{K}\left<x\right> \otimes\mathbb{K}\left<x\right>$ can be written as  
    \begin{equation*}
        c(x^* \otimes x^{k+1}) = \sum_{s+t=k} c_{s,t}  x^s \otimes x^t
    \end{equation*}
    To determine $c_{s,t}$, we use the definition of 1-cocycle in \eqref{cocycle}. \\
    We have in the LHS  
   \begin{equation} \label{LHS n=1}
       c([x^* \otimes x^{k+1}, x^* \otimes x^{l+1}])= c((l-k)x^* \otimes x^{k+l+1})
       =(l-k) \sum_{s+t=l+k} c_{s,t} x^s \otimes x^t
   \end{equation}
   In the RHS, we find 
    \begin{equation*} 
    \begin{aligned}
        &x^* \otimes x^{k+1}(c(x^* \otimes x^{l+1})) - x^* \otimes x^{l+1}(c(x^* \otimes x^{k+1})) \\
        &=x^* \otimes x^{k+1}(\sum_{s+t=l} c_{s,t}  x^s \otimes x^t) - x^* \otimes x^{l+1}(\sum_{s+t=k} c_{s,t}  x^s \otimes x^t) \\
        &=\sum_{s+t=l} c_{s,t} (s x^{s+k} \otimes x^t + t x^{s} \otimes x^{t+k}) - \sum_{s+t=k} c_{s,t} (s x^{s+l} \otimes x^t + t x^{s} \otimes x^{t+l})
    \end{aligned}
    \end{equation*}
    By setting $s \mapsto s-k $ ans $t \mapsto t-k$ in the first sum and $s \mapsto s-l $ ans $t \mapsto t-l$ in the second sum , we obtain 
    \begin{equation} \label{RHS n=1}
    \begin{aligned}
        &\sum_{s+t=l+k} c_{s-k,t} (s-k) x^{s} \otimes x^t + \sum_{s+t=l+k} c_{s,t-k}(t-k) x^{s} \otimes x^{t} - \sum_{s+t=k+l} c_{s-l,t} (s-l) x^{s} \otimes x^t + \sum_{s+t=k+l} c_{s,t-l} (t-l) x^{s} \otimes x^{t}
    \end{aligned}
    \end{equation}
    With \eqref{LHS n=1} and \eqref{RHS n=1}, we find 
    \begin{equation} \label{l-k}
    \begin{aligned}
      &(l-k)c_{s,t}=c_{s-k,t} (s-k) +c_{s,t-k}(t-k)-c_{s-l,t} (s-l) - c_{s,t-l} (t-l) \\
        &\forall l,k \in \mathbb{Z}_{\ge -1} \forall s,t \in \mathbb{Z}_{\ge 0} \text{ s.t. }s+t=l+k.
    \end{aligned}
    \end{equation}
    Note that $c_{s,t}=0 \quad \forall s,t \le 0$. Considering $l>k>0$, we deduce that the coefficients $c_{s,t}$ are determined by those with smaller indices. Thus, the 1-cocycles are fully determined by $c_{0,0},c_{0,1},c_{1,0},c_{1,1},c_{0,2},c_{2,0}$. \\ 
    We can also notice that 
    \begin{equation*}
        \begin{aligned}
            &c(2x^* \otimes x)= 2c_{0,0} 1 \otimes 1\\
            &c([x^*,x^* \otimes x])= c_{0,1} 1\otimes 1 + c_{1,0} 1 \otimes 1
        \end{aligned}
    \end{equation*}
    and
    \begin{equation*}
        \begin{aligned}
            &c(3x^* \otimes x^2)= 3c_{0,1} 1 \otimes x + 3 c_{1,0} x \otimes 1\\
            &c([x^*, x^* \otimes x^3])= x^*(c(x^* \otimes x^3))=2c_{0,2} 1 \otimes x + c_{0,2} x \otimes 1 +c_{1,1} x \otimes 1 + c_{1,1} 1 \otimes 1 .
        \end{aligned}
    \end{equation*}
    This implies that
    \begin{equation} \label{c00}
    \left\{
    \begin{aligned}
    &2c_{0,0}=c_{1,0}+c_{0,1}\\
    &3c_{1,0}=c_{1,1}+2c_{2,0} \\
    &3c_{0,1}= c_{1,1} +2c_{0,2}.
    \end{aligned}
    \right.
    \end{equation}
    Therefore, all coefficients are determined by $c_{1,1},c_{0,2},c_{2,0}$. We now have an upper bound for the number of cocycle which is three.
    The remainder of the proof consists in explicitly describing these three independent cocycles  
    \begin{equation*}
    \left\{
    \begin{aligned}
    &\operatorname{div }\otimes 1: x^* \otimes x^{k+1} \mapsto (k+1)  x^{k} \otimes 1\\ 
    &1 \otimes\operatorname{div } : x^* \otimes x^{k+1} \mapsto (k+1)1 \otimes x^{k}\\
    &\operatorname{Div} : x^* \otimes x^{k+1} \mapsto\sum_{s+t=k}x^s \otimes x^t.
    \end{aligned}
    \right.
    \end{equation*}
    Independence of those three cocycles is clear from coefficients $c_{s,t}$ for s+t=2. But we can observe that 
    \begin{enumerate}
        \item $ \operatorname{div} \otimes 1 $ corresponds to \( c_{s,0} = s+1 \) and \( c_{s,t} = 0 \) if \( t \geq 1 \).
        \item $ 1 \otimes \operatorname{div} $ corresponds to \( c_{0,t} = t+1 \) and \( c_{s,t} = 0 \) if \( s \geq 1 \).
        \item Div corresponds to $c_{s,t}=1 \quad \forall s,t \in \mathbb{Z}_{\geq 0}$.
    \end{enumerate}
\end{Proof*}

\begin{remark*}
    Notice that the equation \eqref{c00} corresponds to choosing  
    \begin{enumerate}
        \item $k=-1, l=1, s=0, t=0$
        \item $k=-1, l=2, s=0, t=1$
        \item $k=-1, l=2, s=1, t=0$
    \end{enumerate}
    in  \eqref{l-k}. 
\end{remark*}

\section{Symplectic derivation Lie algebra of the free algebra} \label{section5} 

In this section, we are focusing on a Lie subalgebra of \(\operatorname{Der}(\mathfrak{L}_{2n}) : \operatorname{Der}_{Sp}(\mathfrak{L}_{2n})\) the symplectic derivation Lie algebra of the free algebra. 
We want to show the uniqueness of \( c: \operatorname{Der}_{Sp}(\mathfrak{L}_{2n}) \to \rvert T(H) \rvert \) a 1-cocycles of degree zero, up to scalar multiple. 
This particular case is of interest due to its similarity in style, albeit in a simpler form. Moreover, it holds significance from the perspective of representation theory.

\subsection{Definitions}

\subsubsection{Free Lie algebra and important results}
Let $\mathfrak{L}_{n}$ denote the free Lie algebra generated by $\{x_1,....,x_n\}$ and $V=\operatorname{Vect}_\mathbb{K}(x_1,...,x_n)$. We can denote \( T(V) =\bigoplus_{i \geq 1}^{\infty} V^{\otimes i} \)the tensor algebra of V. \\
As before, $\mathfrak{L}_{n}$ is naturally graded and has the direct sum decomposition 
\[\mathfrak{L}_{n}= \bigoplus_{k \geq 0}^{\infty} \mathfrak{L}_{n}(k). \] 

\begin{definition}
    A \textbf{derivation} of $\mathfrak{L}_{n}$ is an endomorphisme D of $\mathfrak{L}_{n}$ such that  
\[
D([x_{i_1},x_{i_2}])= [D(x_{i_1}),x_{i_2}] + [x_{i_1},D(x_{i_2})] \quad \forall x_{i_1},x_{i_2} \in \mathfrak{L}_{n}
\]
$\operatorname{Der}(\mathfrak{L}_{n})$ denote the set of all derivation of $\mathfrak{L}_{n}$.
\end{definition}

There exists a direct sum decomposition
\(\operatorname{Der}(\mathfrak{L}_{n})= \bigoplus_{k \geq 0}^{\infty} \operatorname{Der}(\mathfrak{L}_{n})(k) \).

In \cite{Morita}, it was shown that $D \in \operatorname{Der}(\mathfrak{L}_{n})$ is generated by elements of the form \[x^*_i \otimes [[ \cdots [[x_{i_1},x_{i_2}],x_{i_3}], \cdots ]x_{i_{k+1}}] 
\in \operatorname{Der}(\mathfrak{L}_{n})(k).\]
 \begin{remark*}
    Applying the proof method of Theorem \ref{thm:mainresult} to show that there is a unique 1-cocycle of degree zero does not work in this setting, as derivations of degree -1 are not free Lie algebra derivatives.
 \end{remark*}

\begin{definition}
    Let $\mathfrak{L}_{2n}$ be the free Lie algebra generated by $\{x_1,...,x_n,y_1,...,y_n\}$ and $H=\operatorname{Vect} _ \mathbb{K}(x_1,...,x_n,y_1,..., y_n)$. \textbf{The symplectic derivation of Lie algebra of the free algebra} is 
    \begin{equation*}
        \operatorname{Der}_{Sp}(\mathfrak{L}_{2n})=\{ D \in \operatorname{Der}(\mathfrak{L}_{2n}) ; D(\sum_{j=1}^n [x_j,y_j])=0\}.
    \end{equation*}
\end{definition}

\begin{remark*}
    The condition $D\left(\sum_{j=1}^n [x_j, y_j]\right) = 0$ ensures that the derivation $D$ preserves the standard symplectic structure defined on the generating space. 
\end{remark*}

 Hain’s theorem \cite{Hain} shows that, when  $n \to \infty$, $\operatorname{Der}_{Sp}(\mathfrak{L}_{2n})$ is generated by $\wedge^3H$ \cite{ES2} \label{Hain}.
 More precisely, Hain's theorem states that the degree 1 part of the Lie algebra \(\operatorname{Der}_{Sp}(\mathfrak{L}_{2n})\) identifies with \(\wedge^3 H\), and that, in the stable range \(n \to \infty\), the whole Lie algebra \(\operatorname{Der}_{Sp}(\mathfrak{L}_{2n})\) is generated by \(\wedge^3 H\). 

\begin{remark*}
    $\wedge^3H$ is seen as an element of  $\operatorname{Der}_{Sp}(\mathfrak{L}_{2n})$ with the following injection  
    \[
\phi : \wedge^3 H \hookrightarrow \operatorname{Der}_{Sp}(\mathfrak{L}_{2n}) \text{ such that \( \forall z_1,z_2,z_3 \in H\)  }
\]
\begin{equation*}
    \begin{aligned}
        \phi(z_1 \wedge z_2 \wedge z_3) =
\sum_{i=1}^n& x_i^* \otimes \left( \omega(x_i, z_1)[z_2, z_3] + \omega(x_i, z_2)[z_3, z_1] + \omega(x_i, z_3)[z_1, z_2] \right) \\
&+ y_i^* \otimes \left( \omega(y_i, z_1)[z_2, z_3] + \omega(y_i, z_2)[z_3, z_1] + \omega(y_i, z_3)[z_1, z_2] \right)
    \end{aligned}
\end{equation*}

where \(w \) is the symplectic form.
\end{remark*}
We will also use a well-known result in representation theory. 

\begin{theorem} {(Schur's Lemma)} \label{Schur}
     If V and W are irreducible representations of G and $\varphi : V \to W$ is a G-module homomorphism, then 
     \begin{enumerate}
         \item Either $\varphi$ is an isomorphism, or $\varphi =0$. 
         \item If $V=W$, then $\varphi = \lambda Id$ for some $\lambda \in \mathbb{C}$.
     \end{enumerate}
 \end{theorem}

 \begin{remark*}
     For the proof of the theorem and more details, we refer you to \cite{Fulton}. 
 \end{remark*}

\subsubsection{Enomoto-Satoh trace}
The well-known 1-cocycle of degree zero is the Enomoto-Satoh trace \cite{ES1}. To defined it, we we need the following functions \(\forall k \ge 1\)  

\begin{enumerate}
    \item The contraction maps  
        \begin{equation*}
            \begin{aligned}
                \varphi_k : &H^* \otimes H^{\otimes k+1} \to H^{\otimes k} \\
                &z_{0}^* \otimes z_{1} \otimes ... \otimes z_{{k+1}} \mapsto z_{0}^*(z_{1})z_{2} \otimes ... \otimes z_{{k+1}}.
            \end{aligned}
        \end{equation*}
    \item The natural embedding $i_{k+1}$ is the natural inclusion that sends a Lie bracket into the tensor algebra by expanding it recursively  
        \begin{equation*}
            \begin{aligned}
                i_{2} : & \mathfrak{L}_{n}(2) \to H^{\otimes 2} \\
                & [z_{0},z_{1}] \mapsto z_{0} \otimes z_{1} - z_{1} \otimes z_{0}. \\
            \end{aligned}
        \end{equation*}
    \item The projection 
        \begin{equation*}
            \begin{aligned}
             p_k : &H^{\otimes k} \to \rvert T(H) \rvert \\
            & z_{1} \otimes ... \otimes z_{{k}}
            \mapsto      \rvert  z_{1}...z_{{k+1}} \rvert.     
            \end{aligned}
        \end{equation*}
\end{enumerate}

\begin{definition}
    For all k $\ge 1$ the \textbf{contraction map} $ \phi_k$ is a $Gl_n \mathbb{K}$-equivariant homomorphism and is defined as follows 
    \begin{equation*}
        \begin{alignedat}{2}
        \phi_k :\quad & \operatorname{Hom}(H, \mathfrak{L}_{n}(k+1)) := H^* \otimes \mathfrak{L}_{n}(k+1) \to H^{\otimes k} \\
                     & z_{0}^* \otimes [[\cdots[z_{1}, z_{2}], \cdots, z_{{k+1}}] \mapsto
                      \varphi_k \circ (\mathrm{id}_{H^*} \otimes i_{k+1})(z_{0}^* \otimes [[\cdots[z_{1}, z_{2}], \cdots, z_{{k+1}}]). \\
        \end{alignedat}
    \end{equation*}
\end{definition}

\begin{definition} 
    For all k $\ge 2$ the \textbf{Enomoto-Satoh trace} is  
    \begin{equation} \label{def:ES}
        \begin{aligned}
            \operatorname{Tr}_{ES}:= p_k \circ\phi_k : &\operatorname{Hom}(H, \mathfrak{L}_{n}(k+1)) \to \rvert T(H)\rvert.\\
           \end{aligned} \end{equation}
\end{definition}

\begin{remark*}
The Enomoto–Satoh trace is a refinement of Morita's trace 
\[
Tr_M^k  :=
H^* \otimes \mathfrak{L}_{k+1}
\;\xrightarrow\;
S^k(H)
\]
which is a $GL(n,Z)$-equivariant surjective homomorphism. For more information, see \cite{Morita1}, \cite{Morita2} and \cite{ES1}.

\end{remark*}

\subsection{Proof of uniqueness}

In this section, we will demonstrate the theorem \ref{thm:UniqueES}. 

\begin{Proof*}
By \ref{Hain}, the generators of $\operatorname{Der}_{Sp}(\mathfrak{L_{2n}})$ are \(\wedge^3H\) when $ n \to \infty $. We are therefore interested in $\operatorname{dim}(\operatorname{Hom}_{Sp(2n,\mathbb{C})}(\wedge ^3H, T(H)))$ . 

In Fulton and Harris \cite{Fulton}, it is shown that $ \operatorname{Ker} (\bar{\varphi_k})$ is a irreductible representation, where 
  \begin{equation*}
        \begin{aligned}
            \bar{\varphi_k} : &\wedge^k H \to \wedge^{k-2}H\\
            &z_1 \wedge ...\wedge z_k \mapsto \sum_{i < j} w(z_i, z_j) (-1)^{i+j+1}z_1 \wedge ... \wedge \hat{z_i} \wedge ... \wedge \hat{z_j} \wedge ... \wedge z_k 
        \end{aligned}
    \end{equation*}
with $w(z_i, z_j)$ the symplectic form. Moreover, $\bar{\varphi_3}$ is a surjective map 
\begin{equation*}
\begin{aligned}
     \bar{\varphi_3}(x_i \wedge y_i \wedge x_j)&= w(x_i,y_i)x_j - w(x_i,x_j)y_i + w(y_i,x_j)y_i \\
     &= x_j .
\end{aligned}
\end{equation*}

In other words, $\wedge^3 H \to H \to 0$ is an exact split sequence, and we have $\wedge ^3 H \cong H \oplus \operatorname{Ker}(\bar{\varphi_3})$. 

With the Schur's lemma \ref{Schur}, we have $\operatorname{dim}(\operatorname{Hom}_{Sp(2n,\mathbb{C})}(\wedge^3H, H))=1$ : there is a unique homomorphism $Sp$-equivariant from $\wedge ^3H$ to $ H $. As $\operatorname{Der}_{Sp}(\mathfrak{L_{2n}})$ is generated by $\wedge ^3 H$ for $n \to \infty $, we conclude that there exists a unique 1-cocycle \(c: \operatorname{Der}_{Sp}(\mathfrak{L}_{2n}) \to \rvert T(H)\rvert \) of degree zero. 

\end{Proof*}

\section{Appendix} \label{appendix}

For finding the values for $c(x^*_{i_o} \otimes x_{i_1} \otimes x_{i_2} \otimes x_{i_3} \otimes x_{i_4})$, we use the equalities of the proof of the theorem \ref{thm:MoritaSakasai} found in the following paper \cite{Morita}. It gives us for $i_1 \neq i_2 \neq i_3 \neq i_4$ 

\begin{enumerate}
    \item
    \begin{align*}
        c(x^*_{i_o} \otimes x_{i_0} \otimes x_{i_2} \otimes x_{i_3} \otimes x_{i_4}) &= c([x_{i_2}^* \otimes x_{i_2} \otimes x_{i_3} \otimes x_{i_4}, x_{i_0}^* \otimes x_{i_0} \otimes x_{i_2}]) \\
        &= x_{i_2}^* \otimes x_{i_2} \otimes x_{i_3} \otimes x_{i_4}(c(x_{i_0}^* \otimes x_{i_0} \otimes x_{i_2}))- x_{i_0}^* \otimes x_{i_0} \otimes x_{i_2}(c(x_{i_2}^* \otimes x_{i_2} \otimes x_{i_3} \otimes x_{i_4})) \\
        &= x_{i_2}^* \otimes x_{i_2} \otimes x_{i_3} \otimes x_{i_4}( \alpha \rvert x_{i_2} \rvert \otimes 1 + \beta 1 \otimes \rvert x_{i_2} \rvert) - 0 \\
        &= \alpha \rvert x_{i_2}x_{i_3}x_{i_4} \rvert \otimes 1 + \beta 1 \otimes \rvert x_{i_2}x_{i_3}x_{i_4} \rvert.
    \end{align*}
    \item 
    \begin{align*}
         c(x^*_{i_o} \otimes x_{i_1} \otimes x_{i_0} \otimes x_{i_3} \otimes x_{i_4}) &= c([x_{i_0}^* \otimes x_{i_2} \otimes x_{i_3}, x_{i_0}^* \otimes x_{i_1} \otimes x_{i_0} \otimes x_{i_0}]) - c(x_{i_0}^* \otimes x_{i_0} \otimes x_{i_2} \otimes x_{i_3} \otimes x_{i_0}) \\
         &= x_{i_0}^* \otimes x_{i_2} \otimes x_{i_3}(b_1 \rvert x_{i_1}x_{i_0} \rvert \otimes 1 + b_2 1 \otimes \rvert x_{i_1}x_{i_0} \rvert + b_3 \rvert x_{i_1} \rvert \otimes \rvert x_{i_0} \rvert + b_4 \rvert x_{i_0} \rvert \otimes \rvert x_{i_1} \rvert \\
        & + \gamma \rvert x_{i_1}x_{i_0} \rvert \otimes 1 
         + \omega 1 \otimes \rvert x_{i_1}x_{i_0} \rvert) 
         - \gamma \rvert x_{i_1}x_{i_0}x_{i_3} \rvert \otimes 1 - \omega 1 \otimes \rvert x_{i_1}x_{i_0}x_{i_3} \rvert \\
         &=b_1 \rvert x_{i_1}x_{i_3}x_{i_4} \rvert \otimes 1 + b_2 1 \otimes \rvert x_{i_1}x_{i_3}x_{i_4} \rvert + b_3 \rvert x_{i_1} \rvert \otimes \rvert x_{i_3}x_{i_4} \rvert + b_4 \rvert x_{i_3}x_{i_4} \rvert \otimes \rvert x_{i_1} \rvert.
    \end{align*}
 
    \item 
    \begin{align*}
         c(x^*_{i_o} \otimes x_{i_1} \otimes x_{i_2} \otimes x_{i_0} \otimes x_{i_3}) &= c([x_{i_0}^* \otimes x_{i_1} \otimes x_{i_2}, x_{i_0}^* \otimes x_{i_0} \otimes x_{i_0} \otimes x_{i_3}] - x_{i_0}^* \otimes x_{i_0} \otimes x_{i_1} \otimes x_{i_2} \otimes x_{i_3}) \\
         &= c( x_{i_o}^* \otimes x_{i_1} \otimes x_{i_2} ( \alpha \rvert x_{i_0}x_{i_3} \rvert \otimes 1 + \beta 1 \otimes \rvert x_{i_0}x_{i_3} \rvert + b_1 \rvert x_{i_0}x_{i_3} \rvert \otimes 1 + b_2 1 \otimes \rvert   x_{i_0}x_{i_3} \rvert \\
         &+ b_3 \rvert x_{i_0} \rvert \otimes \rvert x_{i_3} \rvert + b_4 \rvert x_{i_3} \rvert \otimes \rvert x_{i_0} \rvert -0) - \alpha \rvert x_{i_1}x_{i_2}x_{i_3} \rvert \otimes 1 - \beta 1 \otimes \rvert x_{i_1}x_{i_2}x_{i_3} \rvert \\
         &=b_1 \rvert x_{i_1}x_{i_2}x_{i_4} \rvert \otimes 1 +b_2 \rvert 1 \otimes \rvert x_{i_1}x_{i_2}x_{i_4} \rvert + b_3 \rvert x_{i_1}x_{i_2} \rvert \otimes \rvert x_{i_4} \rvert +b_4\rvert x_{i_4} \rvert \otimes \rvert x_{i_1}x_{i_2} \rvert.
    \end{align*}
    \item 
    \begin{align*}
        c(x^*_{i_o} \otimes x_{i_1} \otimes x_{i_2} \otimes x_{i_3} \otimes x_{i_0}) &= c([x_{i_0}^* \otimes x_{i_1} \otimes x_{i_2} \otimes x_{i_3}, x_{i_0}^* \otimes x_{i_0} \otimes x_{i_0}]) - c(x^*_{i_0} \otimes x_{i_0} \otimes x_{i_1} \otimes x_{i_2} \otimes x_{i_3}) \\
        &= x_{i_0}^* \otimes x_{i_1} \otimes x_{i_2} \otimes x_{i_3}(\alpha \rvert x_{i_0} \rvert \otimes 1 + \beta \rvert x_{i_0} \rvert \otimes 1 + \gamma \rvert x_{i_0} \rvert \otimes 1 + \omega 1 \otimes \rvert x_{i_0} \rvert ) \\
        &- \alpha \rvert x_{i_1}x_{i_2}x_{i_3} \rvert \otimes 1 + \beta 1 \otimes \rvert x_{i_1}x_{i_2}x_{i_3} \rvert \\
        &= \gamma \rvert x_{i_1}x_{i_2}x_{i_3} \rvert \otimes 1 + \omega 1 \otimes \rvert x_{i_1}x_{i_2}x_{i_3} \rvert.
    \end{align*}

\end{enumerate}

We therefore have as a final formula 
\begin{align*}
    c(x^*_{i_o} \otimes x_{i_1} \otimes x_{i_2} \otimes x_{i_3} \otimes x_{i_4}) &= \delta_{i_0,i_1} (\alpha \rvert x_{i_2}x_{i_3}x_{i_4} \rvert \otimes 1 + \beta 1 \otimes \rvert x_{i_2}x_{i_3}x_{i_4} \rvert) \\
    &+ \delta_{i_0,i_2}(b_1 \rvert x_{i_1}x_{i_3}x_{i_4} \rvert \otimes 1 + b_2 1 \otimes \rvert x_{i_1}x_{i_3}x_{i_4} \rvert + b_3 \rvert x_{i_1} \rvert \otimes \rvert x_{i_3}x_{i_4} \rvert + b_4 \rvert x_{i_3}x_{i_4} \rvert \otimes \rvert x_{i_1} \rvert) \\
    &+\delta_{i_0,i_3}(b_1 \rvert x_{i_1}x_{i_2}x_{i_4} \rvert \otimes 1 +b_2 \rvert 1 \otimes \rvert x_{i_1}x_{i_2}x_{i_4} \rvert + b_3 \rvert x_{i_1}x_{i_2} \rvert \otimes \rvert x_{i_4} \rvert +b_4\rvert x_{i_4} \rvert \otimes \rvert x_{i_1}x_{i_2} \rvert) \\
    &+ \delta_{i_0,i_4} (\gamma \rvert x_{i_1}x_{i_2}x_{i_3} \rvert \otimes 1 + \omega 1 \otimes \rvert x_{i_1}x_{i_2}x_{i_3} \rvert).
\end{align*}

\end{document}